\documentclass{amsart}

\usepackage{amsthm,amsfonts,amsmath,amssymb,mathrsfs,verbatim,cite}

\numberwithin{equation}{section}
\newtheorem{theorem}{Theorem}[section]

\newtheorem{proposition}[theorem]{Proposition}
\newtheorem{lemma}[theorem]{Lemma}

\newcommand{\g}{\mathfrak{g}}

\newcommand{\bil}[2]{(#1,#2)}
\newcommand{\ip}[2]{\langle#1,#2\rangle}
\renewcommand{\Re}{\textup{Re}}
\renewcommand{\and}{\quad\text{and}\quad}
\DeclareMathOperator{\dup}{d\!}
\newcommand{\dupq}{{\dup\,}_q}
\newcommand{\Int}{\int\limits}

\newcommand{\f}{\mathsf{f}}
\newcommand{\Rat}{\mathbb Q}

\newcommand{\Real}{\mathbb R}
\newcommand{\Z}{\mathbb Z}
\newcommand{\F}{\mathbb F}
\newcommand{\Complex}{\mathbb C}
\newcommand{\Symm}{\mathfrak{S}}
\newcommand{\abs}[1]{\lvert#1\rvert}

\newcommand{\la}{\lambda}

\renewcommand{\P}{\mathsf{P}}
\newcommand{\Q}{\mathsf{Q}}
\newcommand{\Pset}{\mathcal{P}}

\newcommand{\sln}[1]{\mathfrak{sl}_{#1}}
\newcommand{\spec}[1]{\langle #1\rangle}

\begin{document}

\title{The $\sln{3}$ Selberg integral}

\author{S. Ole Warnaar}
\thanks{Work supported by the Australian Research Council}

\address{School of Mathematics and Physics, 
The University of Queensland, QLD 4072, Australia}

\date{}

\begin{abstract}
Using an extension of the well-known evaluation symmetry,
a new Cauchy-type identity for Macdonald polynomials is proved.
After taking the classical limit this yields a new $\sln{3}$ 
generalisation of the famous Selberg integral.
Closely related results obtained in this paper are an
$\sln{3}$-analogue of the Askey--Habsieger--Kadell $q$-Selberg integral 
and an extension of the $q$-Selberg integral to 
a transformation between $q$-integrals of different dimensions.
\end{abstract}


\maketitle

\section{Introduction}

Let $\g$ be a simple Lie algebra of rank $n$, with simple roots, 
fundamental weights and Chevalley generators given by $\alpha_i$, 
$\Lambda_i$ and $e_i,f_i,h_i$ for $1\leq i\leq n$.
The roots of $\g$ are normalised such that the maximal root
$\theta$ has length $\sqrt{2}$, i.e., $\bil{\theta}{\theta}=2$,
where $\bil{\cdot\,}{\cdot}$ is the standard bilinear symmetric form on 
the dual of the Cartan subalgebra.

Let $V_{\la}$ and $V_{\mu}$ be highest weight modules of $\g$
with highest weights $\la$ and $\mu$, and denote by
$\text{Sing}_{\la,\mu}[\nu]$ the space of singular 
vectors of weight $\nu$ in $V_{\la}\otimes V_{\mu}$:
\[
\text{Sing}_{\la,\mu}[\nu]=\bigl\{v\in V_{\la}\otimes V_{\mu}:
h_i v=\nu(h_i) v,~e_i v=0,~1\leq i\leq n\bigr\}.
\]

For fixed nonnegative integers $k_1,\dots,k_n$ assign
$k:=k_1+\cdots+k_n$ integration variables $t_1,\dots,t_k$ 
to $\g$ by attaching the $k_i$ variables 
\[
t_{1+k_1+\cdots+k_{i-1}},\dots,t_{k_1+\cdots+k_i}
\]
to the simple root $\alpha_i$. In other words, the first $k_1$ integration
variables are attached to $\alpha_1$, the second $k_2$ to $\alpha_2$ and so on.
By a mild abuse of notation, also set
\[
\alpha_{t_j}=\alpha_i \quad\text{if $k_1+\cdots+k_{i-1}<j\leq k_1+\cdots+k_i$}.
\]

Exploiting the connection between Knizhnik--Zamolodchikov equations
and hypergeometric integrals, see e.g., \cite{EFK03,SV91,Varchenko03},
Mukhin and Varchenko \cite{MV00} conjectured in 2000
that if the space 
\[
\text{Sing}_{\la,\mu}\Bigl[\la+\mu-\sum_{i=1}^n k_i \alpha_i\Bigr]
\]
is one-dimensional, then there exists a real integration domain $\Gamma$ 
such that a closed-form evaluation exists
(in terms of products of ratios of Gamma functions)
for the $\g$ Selberg integral
\begin{equation}\label{Sg}
\Int_{\Gamma}\;
\biggl[\;\:\prod_{i=1}^k t_i^{-\bil{\la}{\alpha_{t_i}}}
(1-t_i)^{-\bil{\mu}{\alpha_{t_i}}}
\prod_{1\leq i<j\leq k}\abs{t_i-t_j}^{\bil{\alpha_{t_i}}{\alpha_{t_j}}}
\biggr]^{\gamma}
\dup t_1\cdots\dup t_k.
\end{equation}

For $\g=\sln{2}$ the evaluation of \eqref{Sg} is well-known and corresponds
to the celebrated Selberg integral \cite{Selberg44,Mehta04,FW08}:
\begin{multline}\label{Selberg}
\Int_{0<t_1<\dots<t_k<1}\,
\prod_{i=1}^k t_i^{\alpha-1}(1-t_i)^{\beta-1}
\prod_{1\le i < j\le k} \abs{t_i-t_j}^{2\gamma}\,
\dup t_1\cdots\dup t_k  \\
=\prod_{i=0}^{k-1} \frac{\Gamma (\alpha+i\gamma)
\Gamma(\beta+i\gamma)\Gamma((i+1)\gamma)}
{\Gamma(\alpha+\beta+(i+k-1)\gamma)\Gamma(\gamma)},
\end{multline}
where
\[
\Re(\alpha)>0,~\Re(\beta)>0,~\Re(\gamma)>
-\min\{1/k,\Re(\alpha)/(k-1),\Re(\beta)/(k-1)\}.
\]
The prospect that generalisations of this extremely important 
integral exist for all simple Lie algebras has led to
much recent progress in evaluating hypergeometric integrals,
see e.g., \cite{FW08,Iguri09,MT05,TV03,Varchenko08,W08,W08b,W09}.

In \cite{TV03} Tarasov and Varchenko obtained an evaluation of
\eqref{Sg} for $\g=\sln{3}$, $\la=\la_1\Lambda_1+\la_2\Lambda_2$, 
$\mu=\mu_2\Lambda_2$ and $k_1\leq k_2$
as follows.
\begin{theorem}[Tarasov--Varchenko]\label{thmTV}
For $0\leq k_1\leq k_2$ let $t=(t_1,\dots,t_{k_1})$,
$s=(s_1,\dots,s_{k_2})$, 
and let $\alpha_1,\alpha_2,\beta_2,\gamma\in\Complex$ such that
$\Re(\alpha_1),\Re(\alpha_2),\Re(\beta_2)>0$ and $\abs{\gamma}$ is 
sufficiently small. Then
\begin{align}\label{STV}
\Int_{C^{k_1,k_2}_{\gamma}[0,1]}
&\prod_{i=1}^{k_1} t_i^{\alpha_1-1}
\prod_{i=1}^{k_2} s_i^{\alpha_2-1} (1-s_i)^{\beta_2-1} \\[-2mm]
&\times \prod_{1\leq i<j\leq k_1} \abs{t_i-t_j}^{2\gamma} 
\prod_{1\leq i<j\leq k_2} \abs{s_i-s_j}^{2\gamma} \:
\prod_{i=1}^{k_1} \prod_{j=1}^{k_2}\, \abs{t_i-s_j}^{-\gamma} \;
\dup t \dup s \notag \\[2mm]
&\qquad=
\prod_{i=0}^{k_1-1}
\frac{\Gamma(\alpha_1+i\gamma)
\Gamma(1+(i-k_2)\gamma)\Gamma((i+1)\gamma)}
{\Gamma(\alpha_1+1+(i+k_1-k_2-1)\gamma) 
\Gamma(\gamma)} \notag \\
&\qquad\quad\times\prod_{i=0}^{k_2-1}
\frac{\Gamma(\alpha_2+i\gamma)
\Gamma(\beta_2+i\gamma)\Gamma((i+1)\gamma)}
{\Gamma(\alpha_2+\beta_2+(i+k_2-k_1-1)\gamma)\Gamma(\gamma)} \notag \\
&\qquad\quad\times\prod_{i=0}^{k_1-1}
\frac{\Gamma(\alpha_1+\alpha_2+(i-1)\gamma)}
{\Gamma(\alpha_1+\alpha_2+\beta_2+(i+k_2-2)\gamma)} \notag \\
&\qquad\quad\times\prod_{i=0}^{k_1-1}
\frac{\Gamma(\alpha_2+\beta_2+(i+k_2-k_1-1)\gamma)} 
{\Gamma(\alpha_2+(i+k_2-k_1)\gamma)}, \notag
\end{align}
where $\dup t=\dup t_1\cdots\dup t_{k_1}$ and
$\dup s=\dup s_1\cdots\dup s_{k_2}$.
\end{theorem}
In the above, $C^{k_1,k_2}_{\gamma}[0,1]$ is a somewhat 
complicated integration chain defined in \eqref{chainTV} on
page~\pageref{chainTV}.
Since
\[
C_{\gamma}^{0,k}[0,1]
=\{(s_1,\dots,s_k)\in\Real^k:~0<s_1<\cdots<s_k<1\}
\]
the Tarasov--Varchenko integral simplifies to the Selberg integral 
when $(k_1,k_2)=(0,k)$.

In \cite{W08,W09} the present author developed a method for
proving Selberg-type integrals using Macdonald polynomials.
This resulted in an evaluation of \eqref{Sg} for
$\g=\sln{n}$ where $\la=\sum_i \la_i\Lambda_i$, $\mu=\mu_n\Lambda_n$
and $k_1\leq k_2\leq \dots\leq k_n$, generalising 
the Selberg and Tarasov--Varchenko integrals.
In this paper we again employ the theory of Macdonald polynomials to establish 
the following Cauchy-type identity. For $\la$ and $\mu$ 
partitions (and not, as above, weights of $\g$) let $\P_{\la}$ be a suitably
normalised Macdonald polynomial. 
Furthermore, let $(a)_n$ be a $q$-shifted factorial 
and $(a)_{\la}$ a generalised $q$-shifted factorial.
(For precise definitions of all of the above, see Section~\ref{SecDN}.)
\begin{theorem}\label{thmCauchy}
Let $X=\{x_1,\dots,x_n\}$ and $Y=\{y_1,\dots,y_m\}$. Then
\begin{multline*}
\sum_{\la,\mu} 
t^{\abs{\la}-n\abs{\mu}}
\P_{\la}(X)\P_{\mu}(Y)\,
(at^{m-1})_{\la}(qt^n/a)_{\mu}
\prod_{i=1}^n \prod_{j=1}^m 
\frac{(a t^{j-i-1})_{\la_i-\mu_j}}
{(a t^{j-i})_{\la_i-\mu_j}} \\
=\prod_{i=1}^n \frac{(ax_i)_{\infty}}{(tx_i)_{\infty}}
\prod_{j=1}^m\frac{(qy_j/a)_{\infty}}{(y_j)_{\infty}}
\prod_{i=1}^n \prod_{j=1}^m \frac{(tx_iy_j)_{\infty}}{(x_iy_j)_{\infty}}.
\end{multline*}
\end{theorem}
For $m=0$ ($n=0$) the above identity reduces to the $q$-binomial
theorem for Macdonald polynomials in $X$ ($Y$).
Theorem~\ref{thmCauchy} may thus be viewed as two coupled, 
multidimensional $q$-binomial theorems. In the special case
$(X,Y,a,q,t)\mapsto (X/q,qY,-q^2,q^2,q^2)$ the theorem simplifies to
Kawanaka's $q$-Cauchy identity for Schur functions \cite{Kawanaka99}
(with the proviso that Kawanaka's description of the summand is 
significantly more involved).

After a limiting procedure, which turns the sums over $\la$ and $\mu$
into integrals, Theorem~\ref{thmCauchy} becomes a new evaluation of 
the Selberg integral \eqref{Sg} for $\g=\sln{3}$ as follows.
\begin{theorem}\label{thmsl3}
Let $t=(t_1,\dots,t_{k_1})$, $s=(s_1,\dots,s_{k_2})$ and let
$\alpha_1,\alpha_2,\beta_1,\beta_2,\gamma\in\Complex$ such that
$\Re(\alpha_1),\Re(\alpha_2),\Re(\beta_1),\Re(\beta_2)>0$,
$\abs{\gamma}$ is sufficiently small, 
\[
\beta_1+(i-k_2-1)\gamma\not\in\Z \quad\text{for $1\leq i\leq \min\{k_1,k_2\}$}
\]
and
\[
\beta_1+\beta_2=\gamma+1.
\]
Then
\begin{align}\label{sl3new}
\Int_{C_{\beta_1,\gamma}^{k_1,k_2}[0,1]}
&\prod_{i=1}^{k_1} t_i^{\alpha_1-1}(1-t_i)^{\beta_1-1}
\prod_{i=1}^{k_2} s_i^{\alpha_2-1} (1-s_i)^{\beta_2-1} \\[-2mm]
&\times \prod_{1\leq i<j\leq k_1} \abs{t_i-t_j}^{2\gamma} 
\prod_{1\leq i<j\leq k_2} \abs{s_i-s_j}^{2\gamma} \:
\prod_{i=1}^{k_1} \prod_{j=1}^{k_2}\, \abs{t_i-s_j}^{-\gamma} \;
\dup t \dup s \notag \\[2mm]
&\qquad=\prod_{i=0}^{k_1-1}\frac{\Gamma(\alpha_1+i\gamma)
\Gamma(\beta_1+(i-k_2)\gamma)\Gamma((i+1)\gamma)}
{\Gamma(\alpha_1+\beta_1+(i+k_1-k_2-1)\gamma)\Gamma(\gamma)} \notag \\
&\qquad\quad\times \prod_{i=0}^{k_2-1}
\frac{\Gamma(\alpha_2+i\gamma)\Gamma(\beta_2+i\gamma)\Gamma((i+1)\gamma)}
{\Gamma(\alpha_2+\beta_2+(i+k_2-k_1-1)\gamma)
\Gamma(\gamma)} \notag \\
&\qquad\quad\times\prod_{i=0}^{k_1-1} 
\frac{\Gamma(\alpha_1+\alpha_2+(i-1)\gamma)}
{\Gamma(\alpha_1+\alpha_2+(i+k_2-1)\gamma)}, \notag 
\end{align}
where $\dup t=\dup t_1\cdots\dup t_{k_1}$,
$\dup s=\dup s_1\cdots\dup s_{k_2}$ and
$C_{\beta,\gamma}^{k_1,k_2}[0,1]$ the integration chain
defined in \eqref{Ckk} on page~\pageref{Ckk}.
\end{theorem}
Since
\[
C_{\beta,\gamma}^{k,0}[0,1]
=\{(t_1,\dots,t_k)\in\Real^k:~0<t_1<\cdots<t_k<1\},
\]
the integral \eqref{sl3new} again contains the Selberg integral 
\eqref{Selberg} as special case.
Unlike \eqref{STV}, however, \eqref{sl3new} exhibits $\Z_2$ symmetry
thanks to
\begin{equation}\label{CCsymm}
C_{\beta_2,\gamma}^{k_2,k_1}[0,1]
=C_{\beta_1,\gamma}^{k_1,k_2}[0,1]
\prod_{i=0}^{k_1-1}\frac{\Gamma(\beta_1+i\gamma)}
{\Gamma(\beta_1+(i-k_2)\gamma)}
\prod_{i=0}^{k_2-1}\frac{\Gamma(\beta_2+(i-k_1)\gamma)}
{\Gamma(\beta_2+i\gamma)}
\end{equation}
for $\beta_1+\beta_2=\gamma+1$,
and
\[
\prod_{i=0}^{k_1-1} 
\frac{\Gamma(\alpha_1+\alpha_2+(i-1)\gamma)}
{\Gamma(\alpha_1+\alpha_2+(i+k_2-1)\gamma)}=
\prod_{i=0}^{k_2-1} 
\frac{\Gamma(\alpha_1+\alpha_2+(i-1)\gamma)}
{\Gamma(\alpha_1+\alpha_2+(i+k_1-1)\gamma)}.
\]

If we specialise $\beta_2=\gamma$ in \eqref{STV}
and $(\beta_1,\beta_2)=(1,\gamma)$ in \eqref{sl3new}
then the respective products over
gamma functions on the right coincide. 
Since also 
\begin{equation}\label{chainid}
C_{1,\gamma}^{k_1,k_2}[0,1]=C_{\gamma}^{k_1,k_2}[0,1]
\end{equation}
(see Section~\ref{Secchains} for more details)
the two $\sln{3}$ integrals are indeed identical for 
this particular specialisation.

\section{Macdonald Polynomials}

\subsection{Definitions and notation}\label{SecDN}

Let $\la=(\la_1,\la_2,\dots)$ be a partition,
i.e., $\la_1\geq \la_2\geq \dots$ with finitely many $\la_i$
unequal to zero.
The length and weight of $\la$, denoted by
$l(\la)$ and $\abs{\la}$, are the number and sum
of the nonzero $\la_i$, respectively.
Two partitions that differ only in their string of zeros
are identified, and the unique partition of length (and weight) 
$0$ is itself denoted by $0$.
The multiplicity of the part $i$ in the partition $\la$ is denoted
by $m_i=m_i(\la)$, and occasionally we will write
$\la=(1^{m_1} 2^{m_2} \dots)$.

We identify a partition with its diagram or Ferrers graph,
defined by the set of points in $(i,j)\in \Z^2$ such that
$1\leq j\leq \la_i$.
The conjugate $\la'$ of $\la$ is the partition obtained by
reflecting the diagram of $\la$ in the main diagonal,
so that, in particular, $m_i(\la)=\la_i'-\la_{i+1}'$.
The statistic $n(\la)$ is given by
\begin{equation*}
n(\la)=\sum_{i\geq 1} (i-1)\la_i=
\sum_{i\geq 1}\binom{\la_i'}{2}.
\end{equation*}

The dominance partial order on the set of partitions of $N$ is
defined by $\la\geq \mu$ if
$\la_1+\cdots+\la_i\geq \mu_1+\cdots+\mu_i$ for all $i\geq 1$.
If $\la\geq \mu$ and $\la\neq\mu$ then $\la>\mu$.

If $\la$ and $\mu$ are partitions then $\mu\subseteq\la$
if (the diagram of) $\mu$ is contained in (the diagram of)
$\la$, i.e., $\mu_i\leq\la_i$ for all $i\geq 1$.

For $s=(i,j)\in\lambda$ the integers
$a(s)$, $a'(s)$, $l(s)$ and $l'(s)$, known as the arm-length, arm-colength,
leg-length and leg-colength of $s$, are defined as
\begin{align*}
a(s)&=\la_i-j,  & a'(s)&=j-1, \\
l(s)&=\la'_j-i,  & l'(s)&=i-1.
\end{align*}
Note that $n(\la)=\sum_{s\in\la}l(s)$.
Using the above we define the generalised hook-length 
polynomials $c_{\la}$ and $c'_{\la}$ as
\begin{align*}
c_{\la}=c_{\la}(q,t)&:=\prod_{s\in\la}\bigl(1-q^{a(s)}t^{l(s)+1}\bigr), \\
c'_{\la}=c'_{\la}(q,t)&:=\prod_{s\in\la}\bigl(1-q^{a(s)+1}t^{l(s)}\bigr).
\end{align*}

The ordinary $q$-shifted factorial are given by
\[
(a)_{\infty}=(a;q)_{\infty}:=\prod_{i=0}^{\infty}(1-aq^i)
\]
and
\[
(b)_z=(a;q)_z:=\frac{(b)_{\infty}}{(bq^z)_{\infty}}.
\]
Note in particular that for $N$ a positive integer 
$(b)_N=(1-b)(1-bq)\cdots(1-bq^{N-1})$, and $1/(q)_{-N}=0$.
Also note that $c'_{(k)}=(q)_k$.
The $q$-shifted factorials can be generalised to allow for a partition
as indexing set: 
\begin{equation*}
(b)_{\la}=(b;q,t)_{\la}:=\prod_{s\in\la}\bigl(1-b q^{a'(s)}t^{-l'(s)}\bigr)
=\prod_{i=1}^{l(\la)}(bt^{1-i})_{\la_i}.
\end{equation*}
With this notation, 
\begin{subequations}\label{ccp}
\begin{align}
c_{\la}&=(t^n)_{\la}
\prod_{1\leq i<j\leq n}\frac{(t^{j-i})_{\la_i-\la_j}}
{(t^{j-i+1})_{\la_i-\la_j}}, \\
c'_{\la}&=(qt^{n-1})_{\la}
\prod_{1\leq i<j\leq n}\frac{(qt^{j-i-1})_{\la_i-\la_j}}
{(qt^{j-i})_{\la_i-\la_j}},
\end{align}
\end{subequations}
where $n$ is an arbitrary integer such that $n\geq l(\la)$.
We also introduce the usual condensed notation 
\begin{equation*}
(a_1,\dots,a_k)_N=(a_1)_N\cdots (a_k)_N
\end{equation*}
and likewise for $q$-shifted factorials indexed by partitions.

\subsection{Macdonald polynomials}
Let $\Symm_n$ denote the symmetric group,
and $\Lambda_n=\Z[x_1,\dots,x_n]^{\Symm_n}$
the ring of symmetric polynomials in $n$ independent variables.

For $X=\{x_1,\dots,x_n\}$ and $\la=(\la_1,\dots,\la_n)$ a partition 
of length at most $n$ the monomial symmetric function $m_{\la}(X)$ 
is defined as
\begin{equation*}
m_{\la}(X)=\sum_{\alpha}x_1^{\alpha_1}\cdots x_n^{\alpha_n},
\end{equation*}
where the sum is over all distinct permutations 
$\alpha=(\alpha_1,\dots,\alpha_n)$ of $\la$.
If $l(\la)>n$ then $m_{\la}(X):=0$.
The monomial symmetric functions $m_{\la}(X)$ for $l(\la)\leq n$
form a $\Z$-basis of $\Lambda_n$.

A $\Rat$-basis of $\Lambda_n$ is given by the power-sum symmetric functions
$p_{\la}(X)$, defined as
\begin{equation*}
p_r(X)=\sum_{i=1}^n x_i^r 
\end{equation*}
for $r\geq 0$ and $p_{\la}(X)=p_{\la_1}(X)\cdots p_{\la_n}(X)$.
The power-sum symmetric functions may be used to define an extremely 
powerful notational tool in symmetric-function theory, known as plethystic 
or $\lambda$-ring notation, see \cite{Haglund08,Lascoux03}.
First we define the plethystic bracket by
\[
f[x_1+\cdots+x_n]=f(x_1,\dots,x_n)
\]
where $f$ is a symmetric function. More simply we just write
\[
f[X]=f(X)
\]
where on the left we assume the additive notation for sets 
(or alphabets), i.e., $X=x_1+\cdots+x_n$ and on the right the 
more conventional $X=\{x_1,\dots,x_n\}$. 
With this notation $f[X+Y]$ takes on the
obvious meaning of the symmetric function $f$ acting on the disjoint union
of the alphabets $X$ and $Y$.
Plethystic notation also allows for the definition of symmetric functions
acting on differences $X-Y$ of alphabets, or for 
symmetric functions acting on such alphabets as $(X-Y)/(1-t)$, 
see e.g., \cite{Lascoux03}.
In this paper we repeatedly need this last alphabet 
when both $X$ and $Y$ contain a single letter, say $a$ and $b$, respectively.
We may then take as definition
\[
p_r\biggl[\frac{a-b}{1-t}\biggr]=\frac{a^r-b^r}{1-t^r},
\]
and extend this by linearity to any symmetric function.
Note in particular that
\[
f\biggl[\frac{1-t^n}{1-t}\biggr]=f(t^{n-1},\dots,t,1)=:f(\spec{0})
\]
corresponds to the so-called principal specialisation, where more generally, 
\[
\spec{\la}=\spec{\la}_n:=
(q^{\la_1}t^{n-1},q^{\la_2}t^{n-2},\dots,q^{\la_n}t^0),
\]
for $l(\la)\leq n$.

\medskip

After this digression we turn to the definition of the Macdonald
polynomials and to some of its basic properties 
\cite{Macdonald88,Macdonald95}.
First we define the scalar product
$\ip{\cdot\,}{\cdot}$ on symmetric functions by
\begin{equation*}
\langle p_{\la},p_{\mu}\rangle=
\delta_{\la\mu} z_{\la} \prod_{i=1}^n
\frac{1-q^{\la_i}}{1-t^{\la_i}},
\end{equation*}
where $z_{\la}=\prod_{i\geq 1} m_i! \: i^{m_i}$ and
$m_i=m_i(\la)$. If we denote the ring of symmetric functions 
in $n$ variables over the field $\F=\Rat(q,t)$ of rational functions 
in $q$ and $t$ by $\Lambda_{n,\F}$, then
the Macdonald polynomial $P_{\la}(X)=P_{\la}(X;q,t)$
is the unique symmetric polynomial in $\Lambda_{n,\F}$ such that:
\[
P_{\la}(X)=m_{\la}(X)+\sum_{\mu<\la}
u_{\la\mu} m_{\mu}(X)
\]
(where $u_{\la\mu}=u_{\la\mu}(q,t)$) and
\begin{equation*}
\langle P_{\la},P_{\mu} \rangle
=0\quad \text{if$\quad\la\neq\mu$.}
\end{equation*}
The Macdonald polynomials $P_{\la}(X)$ with $l(\la)\leq n$ 
form an $\F$-basis of $\Lambda_{n,\F}$. If $l(\la)>n$ then
$P_{\la}(X):=0$.
{}From the definition it follows that $P_{\la}(X)$ for
$l(\la)\leq n$ is homogeneous of (total) degree $\abs{\la}$;
$P_{\la}(zX)=z^{\abs{\la}} P_{\la}(X)$.
A second Macdonald polynomial $Q_{\la}(X)=Q_{\la}(X;q,t)$ is defined
as
\[
Q_{\la}(X)=b_{\la} P_{\la}(X),
\]
where $b_{\la}=b_{\la}(q,t):=c_{\la}/c'_{\la}$. Then
\[
\ip{P_{\la}}{Q_{\mu}}=\delta_{\la\mu}.
\]
This last result may equivalently be stated as the Cauchy identity 
\[
\sum_{\la} P_{\la}(X)Q_{\la}(Y)=
\prod_{i,j=1}^n \frac{(tx_iy_j)_{\infty}}{(x_iy_j)_{\infty}}.
\]

We also need the skew Macdonald polynomials $P_{\la/\mu}(X)$ and 
$Q_{\la/\mu}(X)$ given by
\begin{align*}
P_{\la}[X+Y]&=\sum_{\la}P_{\la/\mu}[X]P_{\mu}[Y] \\
Q_{\la}[X+Y]&=\sum_{\la}Q_{\la/\mu}[X]Q_{\mu}[Y],
\end{align*}
so that $P_{\la/0}(X)=P_{\la}(X)$ and
$Q_{\la/\mu}(X)=b_{\la}b^{-1}_{\mu}P_{\la/\mu}(X)$.
Equivalently,
\[
Q_{\la/\mu}(X)=\sum_{\nu}f_{\mu\nu}^{\la} Q_{\nu}(X),
\]
where $f_{\mu\nu}^{\la}=f_{\mu\nu}^{\la}$ are the 
$q,t$-Littlewood--Richardson coefficients:
\[
P_{\mu}(X)P_{\nu}(X)=\sum_{\la}f_{\mu\nu}^{\la} P_{\la}(X).
\]
{}From the homogeneity of the Macdonald polynomial it
immediately follows that $f_{\mu\nu}^{\la}(q,t)=0$ if 
$\abs{\la}\neq\abs{\mu}+\abs{\nu}$.
It may also be shown that 
$f_{\mu\nu}^{\la}(q,t)=0$ if $\mu,\nu\not\subseteq\la$,
so that $P_{\la/\mu}(X)$ vanishes if $\mu\not\subseteq\la$.

To conclude this section we introduce normalisations of the Macdonald
polynomials convenient for dealing
with basic hypergeometric series with Macdonald polynomial argument:
\begin{subequations}
\begin{align}
\P_{\la/\mu}(X)&=t^{n(\la)-n(\mu)} 
\frac{c'_{\mu}}{c'_{\la}}\, P_{\la/\mu}(X) \\[2mm]
\Q_{\la/\mu}(X)&=t^{n(\mu)-n(\la)} 
\frac{c'_{\la}}{c'_{\mu}}\, Q_{\la/\mu}(X).
\label{skewQP}
\end{align}
\end{subequations}
Note that
\begin{equation}\label{norm}
\Q_{\la/\mu}(X)
=t^{2n(\mu)-2n(\la)} \frac{c_{\la}c'_{\la}}{c_{\mu}c'_{\mu}}\,
\P_{\la/\mu}(X).
\end{equation}
If we also normalise the $q,t$-Littlewood--Richardson coefficients as
\[
\f^{\la}_{\mu\nu}=
t^{n(\mu)+n(\nu)-n(\la)} \frac{c'_{\la}}{c'_{\mu}c'_{\nu}}\, 
f^{\la}_{\mu\nu},
\]
then all of the preceding formulae have perfect analogues:
\begin{subequations}
\begin{align}
\P_{\la}[X+Y]&=\sum_{\mu} \P_{\la/\mu}[Y] \P_{\mu}[X], \\
\Q_{\la}[X+Y]&=\sum_{\mu} \Q_{\la/\mu}[Y] \Q_{\mu}[X],
\label{QQQ}
\end{align}
\end{subequations}
\begin{equation}\label{CXY}
\sum_{\la} \P_{\la}(X)\Q_{\la}(Y)=
\prod_{i,j=1}^n \frac{(tx_iy_j)_{\infty}}{(x_iy_j)_{\infty}},
\end{equation}
\begin{equation}\label{QfQ}
\Q_{\la/\mu}(X)=\sum_{\nu}\f_{\mu\nu}^{\la}\Q_{\nu}(X)
\end{equation}
and
\begin{equation}\label{PPfP}
\P_{\mu}(X)\P_{\nu}(X)=\sum_{\la}\f_{\mu\nu}^{\la}\P_{\la}(X).
\end{equation}

\subsection{Generalised evaluation symmetry}

One of the many striking results in Macdonald polynomial theory ---
first proved in unpublished work by Koornwinder --- 
is the evaluation symmetry
\begin{equation}\label{TK}
\frac{\P_{\la}(\spec{\mu})}{\P_{\la}(\spec{0})}
=\frac{\P_{\mu}(\spec{\la})}{\P_{\mu}(\spec{0})},
\end{equation}
where $\la$ and $\mu$ are partitions of length at most $n$.
As we shall see in Section~\ref{Secproof}, a simple generalisation of this 
result is the key to proving Theorem~\ref{thmCauchy}.
Before stating this generalisation we put \eqref{TK} in plethystic
notation as
\[
\P_{\la}\biggl[\frac{1-t^n}{1-t}\biggr]
\P_{\mu}\bigl[\spec{\la}\bigr]
=\P_{\mu}\biggl[\frac{1-t^n}{1-t}\biggr]
\P_{\la}\bigl[\spec{\mu}\bigr],
\]
where 
\[
f[\spec{\la}]=f\bigl[q^{\la_1}t^{n-1}+\cdots+q^{\la_n}t^0\bigr]=
f\bigl(q^{\la_1}t^{n-1},\dots,q^{\la_n}t^0\bigr)=f(\spec{\la}).
\]
\begin{proposition}[Generalised evaluation symmetry. I]
For $\la$ and $\mu$ partitions of length at most $n$,
\begin{equation}\label{gd}
\P_{\la}\biggl[\frac{1-at^n}{1-t}\biggr]
\P_{\mu}\biggl[a\spec{\la}+\frac{1-a}{1-t}\biggr]
=\P_{\mu}\biggl[\frac{1-at^n}{1-t}\biggr]
\P_{\la}\biggl[a\spec{\mu}+\frac{1-a}{1-t}\biggr].
\end{equation}
\end{proposition}

\begin{proof}
Both sides are polynomials in $a$ of degree 
$\abs{\la}+\abs{\mu}$ with coefficients in $\Rat(q,t)$. 
It thus suffices to verify \eqref{gd} 
for $a=t^p$, where $p$ ranges over the nonnegative integers.
We now write $\spec{\la}=\spec{\la}_n$ and use that
\[
f\biggl[a\spec{\la}_n+\frac{1-a}{1-t}\biggr]\bigg|_{a=t^p}=
f\bigl[\spec{\la}_{n+p}\bigr],
\]
where, since $l(\la)\leq n$,
\[
f\bigl[\spec{\la}_{n+p}\bigr]=
f\bigl(q^{\la_1}t^{n+p-1},\dots,q^{\la_n}t^p,t^{p-1},\dots,t^0\bigr).
\]
As a result we obtain
\[
\P_{\la}\bigl(\spec{0}_{n+p}\bigr) \P_{\mu}\bigl(\spec{\la}_{n+p}\bigr)=
\P_{\mu}\bigl(\spec{0}_{n+p}\bigr) \P_{\la}\bigl(\spec{\mu}_{n+p}\bigr)
\]
which follows from ordinary evaluation symmetry for Macdonald polynomials
on $(n+p)$-letter alphabets.
\end{proof}

The generalised evaluation symmetry can also be stated without resorting
to plethystic notation as a symmetry for skew Macdonald polynomials.
\begin{proposition}[Generalised evaluation symmetry. II]\label{PropII}
\begin{equation}\label{dual}
(at^n)_{\la} \sum_{\nu} (a)_{\nu} \Q_{\mu/\nu}(a\spec{\la})
=(at^n)_{\mu} \sum_{\nu} (a)_{\nu} \Q_{\la/\nu}(a\spec{\mu}).
\end{equation}
\end{proposition}
When $a=1$ both sums vanish unless $\nu=0$. 
Thanks to the principal specialisation formula \cite[page 337]{Macdonald95}
\begin{equation}\label{tn}
\Q_{\la}(\spec{0})=(t^n)_{\la}
\end{equation}
the $a=1$ case of \eqref{dual} thus corresponds to \eqref{TK} in the
equivalent form
\[
\frac{\Q_{\la}(\spec{\mu})}{\Q_{\la}(\spec{0})}
=\frac{\Q_{\mu}(\spec{\la})}{\Q_{\mu}(\spec{0})}.
\]

\begin{proof}[Proof of Proposition~\ref{PropII}]
By changing normalisation we may replace $(\P_{\la},\P_{\mu})$ in
\eqref{gd} by $(\Q_{\la},\Q_{\mu})$.
Using \cite[page 338]{Macdonald95}
\begin{equation}\label{bla}
(a)_{\la}=\Q_{\la}\biggl[\frac{1-a}{1-t}\biggr]
\end{equation}
(for $a=t^n$ this is \eqref{tn}) and \eqref{QQQ}
this gives rise to
\[
(at^n)_{\la}\sum_{\nu}
\Q_{\mu/\nu}\bigl[a\spec{\la}\bigr]
\Q_{\nu}\biggl[\frac{1-a}{1-t}\biggr]
=(at^n)_{\mu}\sum_{\nu} 
\Q_{\la/\nu}\bigl[a\spec{\mu}\bigr]
\Q_{\nu}\biggl[\frac{1-a}{1-t}\biggr].
\]
Once again using \eqref{bla} and dispensing with the remaining plethystic
brackets yields \eqref{dual}.
\end{proof}

\subsection{$\sln{n}$ basic hypergeometric series}

Before we deal with the most important application of the generalised 
evaluation symmetry --- the proof of Theorem~\ref{thmCauchy} ---
we will show how it implies a multivariable generalisation of 
Heine's transformation formula. 
  
Let 
\[
\tau_{\la}=\tau_{\la}(q,t):=(-1)^{\abs{\la}}q^{n(\la')}t^{-n(\la)}
\]
and $X=\{x_1,\dots,x_n\}$. Then the 
$\sln{n}$ basic hypergeometric series $_r\Phi_s$ is defined as
\begin{equation}\label{Phi}
{_r\Phi_s}\biggl[\genfrac{}{}{0pt}{}{a_1,\dots,a_r}
{b_1,\dots,b_s};X\biggr] 
=\sum_{\la} \frac{(a_1,\dots,a_r)_{\la}}{(b_1,\dots,b_s)_{\la}}\,
\tau_{\la}^{s-r+1}\, \P_{\la}(X).
\end{equation}
For $n=1$ this is in accordance with the standard definition
of single-variable basic hypergeometric series ${_r\phi_s}$
as may be found in \cite{AAR99,GR04}:
\begin{align*}
{_r\Phi_s}\biggl[\genfrac{}{}{0pt}{}{a_1,\dots,a_r}
{b_1,\dots,b_s};\{z\}\biggr] 
&=\sum_{k=0}^{\infty} \frac{(a_1,\dots,a_r)_k}{(q,b_1,\dots,b_s)_k}\,
\Bigl((-1)^k q^{\binom{k}{2}}\Bigr)^{s-r+1}\, z^k \\
&={_r\phi_s}\biggl[\genfrac{}{}{0pt}{}{a_1,\dots,a_r}
{b_1,\dots,b_s};z\biggr],
\end{align*}
where in the second line the $q$-dependence of the $_r\phi_s$ series
has been suppressed.

\begin{theorem}[$\sln{n}$--\,$\sln{m}$ transformation formula]
\label{thmPhitrafo}
Let $X=\{x_1,\dots,x_n\}$ and $Y=\{y_1,\dots,y_m\}$. Then
\begin{multline*}
{_{m+1}\Phi_m}\biggl[\genfrac{}{}{0pt}{}{a,ay_1/t,\dots,ay_m/t}
{ay_1,\dots,ay_m};X\biggr] \\
=\biggl(\:\prod_{i=1}^n \frac{(ax_i)_{\infty}}{(x_i)_{\infty}}
\biggr)\biggl(\:\prod_{j=1}^m \frac{(y_j)_{\infty}}{(ay_j)_{\infty}}
\biggr)\,
{_{n+1}\Phi_n}\biggl[\genfrac{}{}{0pt}{}{a,ax_1/t,\dots,ax_n/t}
{ax_1,\dots,ax_n};Y\biggr].
\end{multline*}
\end{theorem}
For $m=0$ this is the $q$-binomial theorem for Macdonald polynomials
\cite{Kaneko96,Macdonald}
\begin{equation}\label{qbt}
{_1\Phi_0}\biggl[\genfrac{}{}{0pt}{}{a}{\text{--}};X\biggr]
=\prod_{i=1}^n \frac{(ax_i)_{\infty}}{(x_i)_{\infty}}
\end{equation}
and for $m=n=1$ it is Heine's $_2\phi_1$ transformation formula
\cite[Equation (III.2)]{GR04}
\[
{_2\phi_1}\biggl[\genfrac{}{}{0pt}{}{a,ay/t}{ay};x\biggr] 
=\frac{(y,ax)_{\infty}}{(x,ay)_{\infty}}\:
{_2\phi_1}\biggl[\genfrac{}{}{0pt}{}{a,ax/t}{ax};y\biggr].
\]

\begin{proof}[Proof of Theorem~\ref{thmPhitrafo}]
First assume that $m=n$, multiply \eqref{dual} by 
$\P_{\la}(X)\P_{\mu}(Y)$ and sum over
$\la$ and $\mu$ to get
\begin{multline}\label{interm}
\sum_{\nu,\mu,\la} (a)_{\nu} (at^n)_{\la} 
\Q_{\mu/\nu}(a\spec{\la}) \P_{\la}(X) \P_{\mu}(Y) \\
=\sum_{\nu,\mu,\la} (a)_{\nu} (at^n)_{\mu}
\Q_{\la/\nu}(a\spec{\mu}) \P_{\la}(X) \P_{\mu}(Y).
\end{multline}
If we multiply \eqref{PPfP} by $\P_{\nu}(Y)$ and sum over $\nu$ then
\eqref{CXY} and \eqref{QfQ} permit this $\nu$-sum to be carried out
explicitly on both sides. As a result we obtain the skew Cauchy identity
(see also \cite[page 352]{Macdonald95})
\begin{equation}\label{Cauchy}
\sum_{\la}\P_{\la}(X)\Q_{\la/\mu}(Y)=
\P_{\mu}(X) \prod_{i,j=1}^n
\frac{(tx_iy_j)_{\infty}}{(x_iy_j)_{\infty}}.
\end{equation}
Applying this to \eqref{interm} we can perform the sum over $\mu$ on the 
left and the sum over $\la$ on the right, leading to
\begin{multline*}
\sum_{\nu,\la} (a)_{\nu} (at^n)_{\la} 
\P_{\la}(X)\P_{\nu}(Y)\prod_{i,j=1}^n 
\frac{(aty_i\spec{\la}_j)_{\infty}}{(ay_i\spec{\la}_j)_{\infty}} \\
=\sum_{\nu,\mu} (a)_{\nu} (at^n)_{\mu}
\P_{\mu}(Y)\P_{\nu}(X)\prod_{i,j=1}^n
\frac{(atx_i\spec{\mu}_j)_{\infty}}{(ax_i\spec{\mu}_j)_{\infty}}.
\end{multline*}
Using the $q$-binomial theorem \eqref{qbt} to perform both sums over
$\nu$ gives
\begin{multline*}
\sum_{\la} (at^n)_{\la} 
\P_{\la}(X) \prod_{i=1}^n \frac{(ay_i)_{\infty}}{(y_i)_{\infty}}
\prod_{i,j=1}^n 
\frac{(aty_i\spec{\la}_j)_{\infty}}{(ay_i\spec{\la}_j)_{\infty}} \\
=\sum_{\mu} (at^n)_{\mu}
\P_{\mu}(Y) \prod_{i=1}^n \frac{(ax_i)_{\infty}}{(x_i)_{\infty}}
\prod_{i,j=1}^n
\frac{(atx_i\spec{\mu}_j)_{\infty}}{(ax_i\spec{\mu}_j)_{\infty}}.
\end{multline*}
Simplifying the products and replacing $a\mapsto at^{-n}$
completes the proof of the theorem for $m=n$.

The general $m,n$ case trivially follows from $m=n$;
assuming without loss of generality
that $m\leq n$ we set $y_{m+1},\dots,y_n=0$
and use that 
\[
\P_{\la}(y_1,\dots,y_m,\underbrace{0,\dots,0}_{n-m})=
\begin{cases}
\P_{\la}(y_1,\dots,y_m)
& \text{if $l(\la)\leq m$}, \\[1mm]
0
& \text{if $l(\la)>m$}.
\end{cases}\qedhere
\]
\end{proof}

\subsection{Proof of Theorem~\ref{thmCauchy}}\label{Secproof}

Using the generalised evaluation symmetry to prove
Theorem~\ref{thmCauchy} is much more difficult than the
proof of Theorem~\ref{thmPhitrafo},
and we proceed by first proving an identity for skew Macdonald
polynomials.

\begin{theorem}\label{thmPQ}
For $\la$ and $\mu$ partitions of length at most $n$,
\begin{multline*}
\sum_{\nu} t^{-\abs{\nu}}
\P_{\mu/\nu}\biggl[\frac{1-a}{1-t}\biggr]
\Q_{\la/\nu}\biggl[\frac{1-q/at}{1-t}\biggr] \\
=t^{-n\abs{\mu}} 
\P_{\mu}\biggl[\frac{1-at^n}{1-t}\biggr]
\Q_{\la}\biggl[\frac{1-qt^{n-1}/a}{1-t}\biggr]
\prod_{i,j=1}^n 
\frac{(qt^{j-i-1}/a)_{\la_i-\mu_j}}{(qt^{j-i}/a)_{\la_i-\mu_j}}.
\end{multline*}
\end{theorem}
Recalling \eqref{bla} and \eqref{norm}
it follows that the right-hand side is completely factorised. 
Moreover, for $a=1$ the summand vanishes unless $\nu=\mu$ so that we
recover the known factorisation of
$\Q_{\la/\mu}[(1-q/t)/(1-t)]$, see \cite[Equation (8.20)]{Rains06} or
\cite[Proposition 3.2]{W05}:
\[
\Q_{\la/\mu}\biggl[\frac{1-q/t}{1-t}\biggr] 
=t^{(1-n)\abs{\mu}} (qt^{n-1})_{\la} \P_{\mu}(\spec{0})
\prod_{i,j=1}^n 
\frac{(qt^{j-i-1})_{\la_i-\mu_j}}{(qt^{j-i})_{\la_i-\mu_j}}.
\]

\begin{proof}
In the first few steps we follow the proof of Theorem~\ref{thmPhitrafo}
but in an asymmetric manner. That is, we take \eqref{dual}, multiply 
both sides by $\P_{\la}(X)$ and sum over $\la$. By the Cauchy
identity \eqref{Cauchy} followed by the $q$-binomial theorem 
\eqref{qbt} we can perform both sums on the right to find
\begin{equation}\label{ass}
\sum_{\la,\nu} (at^n)_{\la} (a)_{\nu} \Q_{\mu/\nu}(a\spec{\la})
\P_{\la}(X)
=(at^n)_{\mu} \prod_{i=1}^n \frac{(ax_i)_{\infty}}{(x_i)_{\infty}}
\prod_{i,j=1}^n \frac{(atx_i\spec{\mu}_j)_{\infty}}
{(ax_i\spec{\mu}_j)_{\infty}}.
\end{equation}
On the left we use \eqref{bla} and \eqref{QQQ} (twice) to rewrite the sum 
over $\nu$ as
\[
\sum_{\nu} (a)_{\nu} \Q_{\mu/\nu}(a\spec{\la})=
\Q_{\mu}\biggl[a\spec{\la}+\frac{1-a}{1-t}\biggr]=
\sum_{\nu} \Q_{\mu/\nu}\biggl[\frac{1-a}{1-t}\biggr]
\Q_{\nu}(a\spec{\la}).
\]
On the right we use \eqref{bla} to trade $(at^n)_{\mu}$ for
$Q_{\mu}[(1-at^n)/(1-t)]$.
Also renaming the summation index $\la$ as $\omega$, \eqref{ass} thus takes
the form
\[
\sum_{\nu,\omega} (at^n)_{\omega} \Q_{\mu/\nu}\biggl[\frac{1-a}{1-t}\biggr]
\Q_{\nu}(a\spec{\omega}) \P_{\omega}(X)
=\Q_{\mu}\biggl[\frac{1-at^n}{1-t}\biggr]
\prod_{i=1}^n \frac{(ax_i)_{\infty}}{(x_i)_{\infty}}
\prod_{i,j=1}^n \frac{(atx_i\spec{\mu}_j)_{\infty}}
{(ax_i\spec{\mu}_j)_{\infty}}.
\]
By \eqref{norm} it readily follows that we may replace all occurrences of
$\Q$ in the above by $\P$. Then specialising $X=b\spec{\la}$ we find
\begin{multline}\label{lhs}
\sum_{\nu,\omega} (at^n)_{\omega} \P_{\mu/\nu}\biggl[\frac{1-a}{1-t}\biggr]
\P_{\nu}(a\spec{\omega}) \P_{\omega}(b\spec{\la}) \\
=\P_{\mu}\biggl[\frac{1-at^n}{1-t}\biggr]
\prod_{i=1}^n \frac{(ab\spec{\la}_i)_{\infty}}
{(b\spec{\la}_i)_{\infty}}
\prod_{i,j=1}^n \frac{(abt\spec{\la}_i\spec{\mu}_j)_{\infty}}
{(ab\spec{\la}_i\spec{\mu}_j)_{\infty}}.
\end{multline}
The next few steps focus on the left-hand side of this identity.
First, by homogeneity followed by an application of the evaluation symmetry 
\eqref{TK},
\[
\text{LHS}\eqref{lhs}=
\sum_{\nu,\omega}a^{\abs{\nu}} b^{\abs{\omega}}(at^n)_{\omega} 
\P_{\mu/\nu}\biggl[\frac{1-a}{1-t}\biggr]
\frac{\P_{\omega}(\spec{0})}{\P_{\la}(\spec{0})}\,
\P_{\la}(\spec{\omega})\P_{\nu}(\spec{\omega}).
\]
Using \eqref{PPfP} this can be further rewritten as
\[
\text{LHS}\eqref{lhs}=
\sum_{\eta,\nu,\omega} a^{\abs{\nu}} b^{\abs{\omega}} 
(at^n)_{\omega} \, \f^{\eta}_{\la\nu}
\P_{\mu/\nu}\biggl[\frac{1-a}{1-t}\biggr]
\frac{\P_{\omega}(\spec{0})}{\P_{\la}(\spec{0})}\,
\P_{\eta}(\spec{\omega}).
\]
By another appeal to evaluation symmetry this yields
\[
\text{LHS}\eqref{lhs}=
\sum_{\eta,\nu,\omega} a^{\abs{\nu}} b^{\abs{\omega}} 
(at^n)_{\omega} \, \f^{\eta}_{\la\nu} 
\P_{\mu/\nu}\biggl[\frac{1-a}{1-t}\biggr]
\frac{\P_{\eta}(\spec{0})}{\P_{\la}(\spec{0})}\,
\P_{\omega}(\spec{\eta}).
\]
The sum over $\omega$ can now be performed by \eqref{qbt} so that
\[
\text{LHS}\eqref{lhs}=
\sum_{\eta,\nu} a^{\abs{\nu}} \, \f^{\eta}_{\la\nu} 
\P_{\mu/\nu}\biggl[\frac{1-a}{1-t}\biggr]
\frac{\P_{\eta}(\spec{0})}{\P_{\la}(\spec{0})}\,
\prod_{i=1}^n \frac{(abt^n\spec{\eta}_i)_{\infty}}
{(b\spec{\eta}_i)_{\infty}}.
\]
Equating this with the right-hand side of \eqref{lhs},  
manipulating the (infinite) $q$-shifted factorials and finally
replacing $b\mapsto bt^{1-n}$ we find
\begin{multline*}
\sum_{\eta,\nu} a^{\abs{\nu}} 
\frac{(b)_{\eta}}{(abt^n)_{\eta}} \,
\f^{\eta}_{\la\nu} 
\P_{\mu/\nu}\biggl[\frac{1-a}{1-t}\biggr]
\P_{\eta}(\spec{0}) \\
=\frac{(b)_{\la}}{(ab)_{\la}} \, \P_{\la}(\spec{0})
\P_{\mu}\biggl[\frac{1-at^n}{1-t}\biggr]
\prod_{i,j=1}^n 
\frac{(abt^{n-i-j+1})_{\la_i+\mu_j}}{(abt^{n-i-j+2})_{\la_i+\mu_j}}.
\end{multline*}
For $a=1$ the summand vanishes unless $\nu=\mu$ and we recover
\cite[Proposition 3.1]{W05}.

Next we specialise $b=q^{-N}$ and then replace $\la$ and $\eta$ by their
complements with respect to the rectangular partition $(N^n)$.
Denoting these complementary partitions by $\hat{\la}$ and $\hat{\eta}$,
we have $\hat{\la}_i=N-\la_{n-i+1}$ for $1\leq i\leq n$ (and a similar
relation between $\eta$ and $\hat{\eta}$).
Using the relations \cite[pp. 259 \& 263]{W05}
\begin{equation*}
\f_{\hat{\la}\nu}^{\hat{\eta}}=
(-q^N t^{1-n})^{\abs{\la-\eta}}
q^{n(\eta')-n(\la')} t^{n(\la)-n(\eta)}
\f_{\eta\nu}^{\la}\,
\frac{(q^{-N})_{\la}}{(q^{-N})_{\eta}}\,
\frac{\P_{\la}(\spec{0})}{\P_{\eta}(\spec{0})},
\end{equation*}
\begin{equation*}
\frac{(a)_{\hat{\la}}}{(b)_{\hat{\la}}}=
\Bigl(\frac{b}{a}\Bigr)^{\abs{\la}}
\frac{(a)_{(N^n)}}{(b)_{(N^n)}}\,
\frac{(q^{1-N}t^{n-1}/b)_{\la}}{(q^{1-N}t^{n-1}/a)_{\la}},
\end{equation*}
and
\[
\P_{\hat{\la}}(\spec{0})=(-1)^{\abs{\la}}
q^{N\abs{\la}-n(\la')}
t^{2N\binom{n}{2}+n(\la)-2(n-1)\abs{\la}} 
\frac{(q^{-N},qt^{n-1})_{\la}}{(qt^{n-1})_{(N^n)}}\,
\P_{\la}(\spec{0}),
\]
as well as the fact that the summand vanishes unless
$\abs{\nu}+\abs{\eta}=\abs{\la}$, we end up with
\begin{multline*}
\sum_{\eta,\nu} t^{-\abs{\nu}}
(q/at)_{\eta}
\f_{\eta\nu}^{\la}\,
\P_{\mu/\nu}\biggl[\frac{1-a}{1-t}\biggr] \\
=t^{-n\abs{\mu}} (qt^{n-1}/a)_{\la}
\P_{\mu}\biggl[\frac{1-at^n}{1-t}\biggr]
\prod_{i,j=1}^n 
\frac{(qt^{j-i-1}/a)_{\la_i-\mu_j}}{(qt^{j-i}/a)_{\la_i-\mu_j}}.
\end{multline*}
By \eqref{QfQ} and \eqref{bla}
\begin{equation}\label{bfQ}
\sum_{\nu} (b)_{\nu} \f^{\la}_{\mu\nu}=
\sum_{\nu} \f^{\la}_{\mu\nu} \Q_{\nu}\biggl[\frac{1-b}{1-t}\biggr]=
\Q_{\la/\mu}\biggl[\frac{1-b}{1-t}\biggr],
\end{equation}
so that the sum over $\eta$ can be performed. By a final appeal to
\eqref{bla} the proof is done.
\end{proof}

Equipped with Theorem~\ref{thmPQ} it is not difficult to prove 
Theorem~\ref{thmCauchy}. To streamline the proof given below
we first prepare an easy lemma.
\begin{lemma}\label{Pieri}
For $X=\{x_1,\dots,x_n\}$ and $\mu$ a partition of length at most $n$,
\begin{align*}
\sum_{\la} \Q_{\la/\mu}\biggl[\frac{a-b}{1-t}\biggr]\P_{\la}(X)
&=\P_{\mu}(X) \prod_{i=1}^n  \frac{(bx_i)_{\infty}}{(ax_i)_{\infty}} \\
\intertext{and} 
\sum_{\la} \P_{\la/\mu}\biggl[\frac{a-b}{1-t}\biggr]\Q_{\la}(X)
&=\Q_{\mu}(X) \prod_{i=1}^n  \frac{(bx_i)_{\infty}}{(ax_i)_{\infty}}.
\end{align*}
\end{lemma}
For $\mu=0$ this is just the $q$-binomial theorem \eqref{qbt}
for Macdonald polynomials.

\begin{proof}
By \eqref{norm} the two identities are in fact one and the same 
result and we only need to prove the first claim.
To achieve this we multiply \eqref{bfQ} by $\P_{\la}(aX)$ and sum over $\la$.
By \eqref{PPfP} and homogeneity this yields 
\[
\P_{\mu}(X)
\sum_{\nu} (b)_{\nu} \P_{\nu}(aX)=
\sum_{\la} \Q_{\la/\mu}\biggl[\frac{a-ab}{1-t}\biggr]\P_{\la}(X).
\]
On the right we can sum over $\nu$ using the $q$-binomial theorem \eqref{qbt}
leading to the desired result (with $b\mapsto ab$).
\end{proof}

\begin{proof}[Proof of Theorem~\ref{thmCauchy}]
Elementary manipulations show that the theorem is invariant 
under the simultaneous changes
$n\leftrightarrow m$, $tX\leftrightarrow Y$ and $a\mapsto qt/a$.
Without loss of generality we may thus assume that $m\leq n$. But
\[
(at^{n-1})_{\la}
\prod_{i,j=1}^n \frac{(a t^{j-i-1})_{\la_i-\mu_j}}
{(a t^{j-i})_{\la_i-\mu_j}}\bigg|_{\mu_{m+1}=\dots=\mu_n=0}
=(at^{m-1})_{\la}
\prod_{i=1}^n \prod_{j=1}^m \frac{(a t^{j-i-1})_{\la_i-\mu_j}}
{(a t^{j-i})_{\la_i-\mu_j}}
\]
so that the case $m<n$ follows from the case $m=n$ by setting
$y_{m+1}=\dots=y_n=0$. 

In the remainder we assume that $m=n$, 
in which case the theorem simplifies to
\begin{multline}\label{misn}
\sum_{\la,\mu} 
t^{\abs{\la}-n\abs{\mu}}
\P_{\la}(X)\P_{\mu}(Y)\,
(at^{n-1})_{\la}(qt^n/a)_{\mu}
\prod_{i,j=1}^n 
\frac{(a t^{j-i-1})_{\la_i-\mu_j}}
{(a t^{j-i})_{\la_i-\mu_j}} \\
=\prod_{i=1}^n \frac{(ax_i)_{\infty}}{(tx_i)_{\infty}}
\prod_{j=1}^n\frac{(qy_j/a)_{\infty}}{(y_j)_{\infty}}
\prod_{i,j=1}^n \frac{(tx_iy_j)_{\infty}}{(x_iy_j)_{\infty}}.
\end{multline}
To prove this we take Theorem~\ref{thmPQ},
replace $a\mapsto q/a$, multiply both sides by
\[
t^{\abs{\la}} \P_{\la}(X)\Q_{\mu}(Y)
\]
and sum over $\la$ and $\mu$. Hence
\begin{multline*}
\sum_{\la,\mu,\nu} \P_{\la}(X)\Q_{\mu}(Y)
\P_{\mu/\nu}\biggl[\frac{1-q/a}{1-t}\biggr]
\Q_{\la/\nu}\biggl[\frac{t-a}{1-t}\biggr] \\
=\sum_{\la,\mu} t^{\abs{\la}-n\abs{\mu}} 
\P_{\la}(X)\Q_{\mu}(Y) (at^{n-1})_{\la}\,
\P_{\mu}\biggl[\frac{1-qt^n/a}{1-t}\biggr]
\prod_{i,j=1}^n 
\frac{(at^{j-i-1})_{\la_i-\mu_j}}{(at^{j-i})_{\la_i-\mu_j}}.
\end{multline*}
On the right we apply \eqref{norm} and \eqref{bla} to rewrite
\[
\Q_{\mu}(Y) \P_{\mu}\biggl[\frac{1-qt^n/a}{1-t}\biggr]=
(qt^n/a)_{\mu} \P_{\mu}(Y),
\]
and on the left we employ Lemma~\ref{Pieri} to carry out the sums 
over $\la$ and $\mu$. Hence
\begin{multline*}
\sum_{\nu} \P_{\nu}(X)\Q_{\nu}(Y)
\prod_{i=1}^n \frac{(ax_i)_{\infty}}{(tx_i)_{\infty}}
\prod_{j=1}^n \frac{(qy_j/a)_{\infty}}{(y_j)_{\infty}} \\
=\sum_{\la,\mu} t^{\abs{\la}-n\abs{\mu}} 
\P_{\la}(X)\P_{\mu}(Y) (at^{n-1})_{\la}(qt^n/a)_{\mu} \,
\prod_{i,j=1}^n 
\frac{(at^{j-i-1})_{\la_i-\mu_j}}{(at^{j-i})_{\la_i-\mu_j}}.
\end{multline*}
Performing the remaining sum on the left by \eqref{CXY} 
results in \eqref{misn}.
\end{proof}

We conclude this section with a remark about a generalisation of
Theorem~\ref{thmCauchy}.
Let $X=\{x_1,\dots,x_n\}$ and let $\la$ be a partition of length $n$.
Then
\[
P_{\la}(X)=x_1\cdots x_n \, P_{\mu}(X),
\]
where $\mu=(\la_1-1,\dots,\la_n-1)$.
Now let $\Pset$ denote the set of weakly decreasing integer sequences of
finite length. Then we may turn things around and
use the above recursion to extend
$P_{\la}$ to all $\la\in\Pset$.
It is then readily verified that
\[
(qt^{n-1})_{\la} \P_{\la}(X)
= t^{n(\la)} \frac{(qt^{n-1})_{\la}}{c'_{\la}}\,
P_{\la}(X)
\]
is well-defined for $\la\in\Pset$ (unlike $\P_{\la}(X)$).

We now state without proof the following generalisation of 
Theorem~\ref{thmCauchy}.
\begin{theorem}
Let $X=\{x_1,\dots,x_n\}$ and $Y=\{y_1,\dots,y_m\}$. Then
\begin{multline*}
\sum_{\la\in\Pset}\sum_{\mu}  t^{\abs{\la}-n\abs{\mu}}
\P_{\la}(X) \P_{\mu}(Y) 
\frac{(at^{m-1},qt^{n-1})_{\la}(bt^n)_{\mu}}{(abt^{n-1})_{\la}}\,
\prod_{i=1}^n \prod_{j=1}^m
\frac{(at^{j-i-1})_{\la_i-\mu_j}}{(at^{j-i})_{\la_i-\mu_j}} \\
=\prod_{i=1}^n\frac{(qt^{i-1},bt^i,ax_i,q/ax_i)_{\infty}}
{(abt^{i-1},qt^i/a,tx_i,b/x_i)_{\infty}}
\prod_{j=1}^m  \frac{(by_j)_{\infty}}{(y_j)_{\infty}} 
\prod_{i=1}^n \prod_{j=1}^m \frac{(tx_iy_j)_{\infty}}{(x_iy_j)_{\infty}}.
\end{multline*}
\end{theorem}
For $m=0$ this reduces to Kaneko's $_1\Psi_1$ sum for Macdonald
polynomials \cite{Kaneko98} and for $n=0$ to the $_1\Phi_0$ sum 
\eqref{qbt}.
When $ab=q$ the summand on the left vanishes unless
$\la$ is an actual partition and we recover Theorem~\ref{thmCauchy}.

\section{The $\sln{3}$ Selberg integral}\label{SecSel3}

\subsection{The integration chains $C^{k_1,k_2}_{\beta,\gamma}[0,1]$ and
$C^{k_1,k_2}_{\gamma}[0,1]$}\label{Secchains}

Before proving Theorem~\ref{thmsl3} 
we give two descriptions of the chain $C_{\beta,\gamma}^{k_1,k_2}[0,1]$.
We also identify the special case $\beta=1$ with the
chain $C_{\gamma}^{k_1,k_2}[0,1]$ defined by Tarasov and Varchenko in 
\cite{TV03}.

Let
\begin{multline}\label{Rk1k2}
I^{k_1,k_2}[0,1]=\{(x_1,\dots,x_{k_1},y_1,\dots,y_{k_2})\in\Real^{k_1+k_2}:\\
0<x_1<\dots<x_{k_1}<1\and
0<y_1<\dots<y_{k_2}<1\},
\end{multline}
and fix a total ordering among the $x_i$ and $y_j$ as follows.
Let $a=(a_1,\dots,a_{k_1})$ be a weakly increasing sequence of nonnegative
integers not exceeding $k_2$:
\begin{equation}\label{aas}
0\leq a_1\leq\dots\leq a_{k_1}\leq k_2.
\end{equation}
Then the domain $I_a^{k_1,k_2}[0,1]\subseteq I^{k_1,k_2}[0,1]$
is formed by imposing the additional inequalities
\[
x_i<y_{a_i+1}<y_{a_i+2}<\dots<y_{a_{i+1}}<x_{i+1} \quad
\text{for $\;0\leq i\leq k_1$},
\]
where $x_0:=0$, $x_{k_1+1}:=1$, $a_0:=0$ and $a_{k_1+1}:=k_2$. 
Equivalently,
\begin{equation}\label{order}
\begin{cases}
0<y_1<y_2<\dots<y_{a_i}<x_i \\
x_i<y_{a_i+1}<\dots<y_{k_2-1}<y_{k_2}<1
\end{cases}
\quad \text{for $\;1\leq i\leq k_1$}.
\end{equation}
Clearly, as a chain,
\begin{equation}\label{Ikk}
I^{k_1,k_2}[0,1]=\sum_a I_a^{k_1,k_2}[0,1],
\end{equation}
where the sum is over all sequences $a=(a_1,\dots,a_{k_1})$
satisfying \eqref{aas}.
To lift $I^{k_1,k_2}[0,1]$ to $C^{k_1,k_2}_{\beta,\gamma}[0,1]$
we replace the right-hand side of \eqref{Ikk} by a weighted sum:
\begin{equation}\label{Ckk}
C^{k_1,k_2}_{\beta,\gamma}[0,1]=\sum_{a}
\biggl(\:\prod_{i=1}^{k_1} 
\frac{\sin\pi(\beta-(i-a_i-k_1+k_2)\gamma)}
{\sin\pi(\beta-(i-k_1+k_2)\gamma)} \biggr) I_a^{k_1,k_2}[0,1],
\end{equation}
where it is assumed that $\beta,\gamma\in\Complex$ such that
\[
\beta+(i-k_2-1)\gamma\not\in\Z\quad 
\text{for $1\leq i\leq \min\{k_1,k_2\}$.}
\]
This is a necessary and sufficient condition
for 
\[
\prod_{i=1}^{k_1} 
\frac{\sin\pi(\beta-(i-a_i-k_1+k_2)\gamma)}
{\sin\pi(\beta-(i-k_1+k_2)\gamma)}
\]
to be free of poles for all admissible sequences $a$.

By viewing $(a_{k_1},\dots,a_2,a_1)$ as
a partition with largest part not exceeding $k_2$ and length
not exceeding $k_1$, the operations of conjugation and/or
complementation yield several alternative descriptions of the 
chain \eqref{Ckk}.
Below we give one such description, reflecting
the $\Z_2$ symmetry of Theorem~\ref{thmsl3} with respect
to the interchange of the labels $1$ and $2$ in $k_i$, $\alpha_i$ and
$\beta_i$. 

Assume \eqref{Rk1k2} and fix a total ordering among the $x_i$ and $y_j$ as
follows.
Let $b=(b_1,\dots,b_{k_2})$ be a weakly increasing sequence of nonnegative
integers not exceeding $k_1$:
\begin{equation}\label{bb}
0\leq b_1\leq\dots\leq b_{k_2}\leq k_1.
\end{equation}
Then the domain $\bar{I}_b^{k_1,k_2}[0,1]\subseteq I^{k_1,k_2}[0,1]$
is formed by assuming the further inequalities 
\[
y_i<x_{b_i+1}<x_{b_i+2}<\dots<x_{b_{i+1}}<y_{i+1} \quad
\text{for $\;0\leq i\leq k_2$},
\]
where $y_0:=0$, $y_{k_2+1}:=1$, $b_0:=0$ and $b_{k_2+1}:=k_1$. 
It is easily seen that if $\mu=(b_{k_2},\dots,b_1)$
and $\la=(a_{k_1},\dots,a_1)$, then $\mu'$ is the conjugate of
$\la$ with respect to $(k_2^{k_1})$, i.e., 
$\mu'_i=k_2-\la_{k_1-i+1}=k_2-a_i$
for $1\leq i\leq k_1$.
Hence, for a pair of admissible sequences $(a,b)$
related by ``conjugation--complementation'',
\[
\bar{I}_b^{k_1,k_2}[0,1]=I_a^{k_1,k_2}[0,1]
\]
and
\[
\prod_{i=1}^{k_2} 
\frac{\sin\pi(\beta+(i-b_i+k_1-k_2-1)\gamma)}
{\sin\pi(\beta+(i-k_2-1)\gamma)}
=
\prod_{i=1}^{k_1} 
\frac{\sin\pi(\beta-(i-a_i-k_1+k_2)\gamma)}
{\sin\pi(\beta-(i-k_1+k_2)\gamma)}.
\]
In other words,
\begin{equation}\label{Ckkb}
C^{k_1,k_2}_{\beta,\gamma}[0,1]
=\sum_b \biggl(\:\prod_{i=1}^{k_2} 
\frac{\sin\pi(\beta+(i-b_i+k_1-k_2-1)\gamma)}
{\sin\pi(\beta+(i-k_2-1)\gamma)} \biggr)
\bar{I}_b^{k_1,k_2}[0,1]
\end{equation}
summed over all sequences $b=(b_1,\dots,b_{k_2})$ subject to
\eqref{bb}.
Comparing \eqref{Ckk} and \eqref{Ckkb}, and using that 
for $\beta_1+\beta_2=\gamma+1$,
\[
\prod_{i=1}^{k_1} 
\frac{\sin\pi(\beta_1-(i-k_1+k_2))\gamma)}
{\sin\pi(\beta_2+(i-k_1-1)\gamma)}
=
\prod_{i=0}^{k_1-1}\frac{\Gamma(\beta_1+i\gamma)}
{\Gamma(\beta_1+(i-k_2)\gamma)}
\prod_{i=0}^{k_2-1}\frac{\Gamma(\beta_2+(i-k_1)\gamma)}
{\Gamma(\beta_2+i\gamma)}
\]
it readily follows that the symmetry relation \eqref{CCsymm} holds.

\medskip

To conclude this section we consider \eqref{Ckk} for $\beta=1$:
\[
C^{k_1,k_2}_{1,\gamma}[0,1]=\sum_a
\biggl(\:\prod_{i=1}^{k_1} 
\frac{\sin\pi((i-a_i-k_1+k_2)\gamma)}
{\sin\pi((i-k_1+k_2)\gamma)} \biggr)
I_a^{k_1,k_2}[0,1].
\]
The summand vanishes if $a_i=i-k_1+k_2$ for some
$1\leq i\leq k_1$ so that we may add the additional
restrictions 
\[
a_i\neq i-k_1+k_2 \qquad\text{for $1\leq i\leq k_1$}
\]
to the sum over $a$.
Recalling \eqref{aas} this in fact implies that 
the much stronger 
\[
a_i<i-k_1+k_2 \qquad\text{for $1\leq i\leq k_1$}.
\]
Therefore,
\[
C^{k_1,k_2}_{1,\gamma}[0,1]=
\sum_{\substack{a \\[1pt] a_i<i-k_1+k_2}}
\biggl(\:\prod_{i=1}^{k_1} 
\frac{\sin\pi((i-a_i-k_1+k_2)\gamma)}
{\sin\pi((i-k_1+k_2)\gamma)} \biggr)
I_a^{k_1,k_2}[0,1].
\]
Defining $M(i)=a_i+1$,
so that
\[
1\leq M(1)\leq M(2)\leq \dots \leq M(k_1)\leq k_2
\]
and
\[
M(i)\leq i-k_i+k_2 \qquad\text{for $1\leq i\leq k_1$},
\]
and writing $M=(M(1),\dots,M(k_1))$, we finally obtain
\begin{align}\label{chainTV}
C^{k_1,k_2}_{1,\gamma}[0,1]
&=\sum_M
\biggl(\:\prod_{i=1}^{k_1} 
\frac{\sin\pi((i-M(i)-k_1+k_2+1)\gamma)}
{\sin\pi((i-k_1+k_2)\gamma)} \biggr)
I_M^{k_1,k_2}[0,1] \\[1mm]
&=:C^{k_1,k_2}_{\gamma}[0,1]. \notag
\end{align}
In the above, by abuse of notation, $I_M^{k_1,k_2}[0,1]=I_a^{k_1,k_2}[0,1]$
if $M=(a_1+1,\dots,a_{k_1}+1)$. The chain $C^{k_1,k_2}_{\gamma}[0,1]$ is
precisely that of Tarasov and Varchenko (up to an
interchange of $k_1$ and $k_2$), see \cite[page 177]{TV03}.

\subsection{Proof of Theorem~\ref{thmsl3}}
We are now prepared to prove Theorem~\ref{thmsl3}. In fact, we will prove
a more general integral, generalising Kadell's extension of the Selberg 
integral \cite{Kadell97} to $\sln{3}$.
To shorten some of the subsequent equations we introduce another
normalised Macdonald polynomial, and for $X=\{x_1,\dots,x_n\}$
\[
\tilde{\P}_{\la}(X)=\frac{\P_{\la}(X)}{\P_{\la}(\spec{0})}=
\frac{P_{\la}(X)}{P_{\la}(\spec{0})}.
\]
Similarly we define a (normalised) Jack polynomial as
\[
\tilde{\P}_{\la}^{(\alpha)}(X)=
\lim_{q\to 1} \tilde{\P}_{\la}(X;q^{\alpha},q).
\]
Hence
\[
\tilde{\P}_{\la}^{(\alpha)}(X)=
\frac{P_{\la}^{(\alpha)}(X)}{P_{\la}^{(\alpha)}(1^n)},
\]
where $P_{\la}^{(\alpha)}(X)$ is the Jack polynomial 
\cite{Macdonald95,Stanley89}.

\begin{theorem}
Set $X=\{x_1,\dots,x_{k_1}\}$, $Y=\{y_1,\dots,y_{k_2}\}$,
\[
\dup X=\dup x_1\cdots\dup x_{k_1}\and
\dup Y=\dup y_1\cdots\dup y_{k_1}.
\]
For $\alpha_1,\alpha_2,\beta_1,\beta_2,\gamma\in\Complex$ such that
$\abs{\gamma}$ is sufficiently small, 
\[
\min\{\Re(\alpha_1)+\la_{k_1},\Re(\alpha_2)+\mu_{k_2},
\Re(\beta_1),\Re(\beta_2)\}>0,
\]

\[
\beta_1+(i-k_2-1)\gamma\not\in\Z \quad\text{for $1\leq i\leq \min\{k_1,k_2\}$}
\]
and
\[
\beta_1+\beta_2=\gamma+1
\]
there holds
\begin{align*}
\Int_{C_{\beta_1,\gamma}^{k_1,k_2}[0,1]}
&\tilde{\P}^{(1/\gamma)}_{\la}(X)
\tilde{\P}^{(1/\gamma)}_{\mu}(Y)
\prod_{i=1}^{k_1} x_i^{\alpha_1-1}(1-x_i)^{\beta_1-1}
\prod_{i=1}^{k_2} y_i^{\alpha_2-1} (1-y_i)^{\beta_2-1} \\[-2mm]
&\times \prod_{1\leq i<j\leq k_1} \abs{x_i-x_j}^{2\gamma} 
\prod_{1\leq i<j\leq k_2} \abs{y_i-y_j}^{2\gamma} \:
\prod_{i=1}^{k_1} \prod_{j=1}^{k_2}\, \abs{x_i-y_j}^{-\gamma} \;
\dup X \dup Y \\[2mm]
&\qquad=\prod_{i=1}^{k_1}\frac{\Gamma(\alpha_1+(k_1-i)\gamma+\la_i)
\Gamma(\beta_1+(i-k_2-1)\gamma)\Gamma(i\gamma)}
{\Gamma(\alpha_1+\beta_1+(2k_1-k_2-i-1)\gamma+\la_i)\Gamma(\gamma)} \\
&\qquad\quad\times \prod_{i=1}^{k_2}
\frac{\Gamma(\alpha_2+(k_2-i)\gamma+\mu_i)
\Gamma(\beta_2+(i-1)\gamma)\Gamma(i\gamma)}
{\Gamma(\alpha_2+\beta_2+(2k_2-k_1-i-1)\gamma+\mu_i)
\Gamma(\gamma)} \\
&\qquad\quad\times\prod_{i=1}^{k_1} \prod_{j=1}^{k_2} 
\frac{\Gamma(\alpha_1+\alpha_2+(k_1+k_2-i-j-1)\gamma+\la_i+\mu_j)}
{\Gamma(\alpha_1+\alpha_2+(k_1+k_2-i-j)\gamma+\la_i+\mu_j)}.
\end{align*}
\end{theorem}
Theorem~\ref{thmsl3} corresponds to the case special case $\la=\mu=0$,
and Kadell's integral arises by taking $k_1=0$ or $k_2=0$.

\begin{proof}
Throughout the proof $0<q<1$.

We take Theorem~\ref{thmCauchy} with $(\la,\mu)$ replaced
by $(\eta,\nu)$ and $(n,m)$ replaced by $(k_1,k_2)$. 
If we then specialise $X=z\spec{\la}_{k_1}$ and $Y=w\spec{\mu}_{k_2}$ 
and use the evaluation symmetry \eqref{TK} on both Macdonald polynomials 
in the summand, we obtain
\begin{multline*}
\sum_{\eta,\nu} 
t^{\abs{\eta}-k_1\abs{\nu}}
\tilde{\P}_{\la}(\spec{\eta}_{k_1})
\tilde{\P}_{\mu}(\spec{\nu}_{k_2})
\P_{\eta}(z\spec{0}_{k_1})\P_{\nu}(w\spec{0}_{k_2}) \\[-2mm]
\times (at^{k_2-1})_{\la}(qt^{k_1}/a)_{\nu}
\prod_{i=1}^{k_1} \prod_{j=1}^{k_2} 
\frac{(a t^{j-i-1})_{\eta_i-\nu_j}}{(a t^{j-i})_{\eta_i-\nu_j}} \\
=\prod_{i=1}^{k_1} \frac{(azq^{\la_i}t^{k_1-i})_{\infty}}
{(zq^{\la_i}t^{k_1-i+1})_{\infty}}
\prod_{j=1}^{k_2}\frac{(wq^{\mu_j+1}t^{k_2-j}/a)_{\infty}}
{(wq^{\mu_j}t^{k_2-j})_{\infty}}
\prod_{i=1}^{k_1} \prod_{j=1}^{k_2} 
\frac{(wztq^{\la_i+\mu_j}t^{k_1+k_2-i-j})_{\infty}}
{(wzq^{\la_i+\mu_j}t^{k_1+k_2-i-j})_{\infty}}.
\end{multline*}
Next we set
\[
(z,w,a,t)=(q^{\alpha_1-\gamma},q^{\alpha_2},
q^{\beta_1+(k_1-k_2)\gamma},q^{\gamma})
\]
and introduce the auxiliary variable $\beta_2$ by
$\beta_1+\beta_2=\gamma+1$.
Equations \eqref{ccp}, \eqref{norm} and \eqref{bla} 
imply the principal specialisation formula
\[
\P_{\la}(\spec{0}_n)=\frac{t^{2n(\la)}}{(qt^{n-1})_{\la}}
\prod_{1\leq i<j\leq n}
\frac{1-q^{\la_i-\la_j}t^{j-i}}{1-t^{j-i}}\,
\frac{(t^{j-i+1})_{\la_i-\la_j}}{(qt^{j-i-1})_{\la_i-\la_j}}.
\]
Using this as well as the definition of the $q$-Gamma function
\[
\Gamma_q(x)=\frac{(q)_{x-1}}{(1-q)^{x-1}}, \quad x\in\Complex,
\]
we can rewrite the above identity as
\begin{multline*}
(1-q)^{k_1+k_2}\sum_{\eta,\nu} 
\tilde{\P}_{\la}(x_1q^{(k_1-1)\gamma},x_2q^{(k_1-2)\gamma},\dots,x_{k_1}) \\
\times
\tilde{\P}_{\mu}(y_1q^{(k_2-1)\gamma},y_2q^{(k_2-2)\gamma},\dots,y_{k_1}) \\
\times
\prod_{i=1}^{k_1} x_i^{\alpha_1}
(q^{1+(k_1-i)\gamma}x_i)_{\beta_1-1} 
\prod_{1\leq i<j\leq k_1} x_j^{2\gamma}
\bigl(1-q^{(j-i)\gamma}x_i/x_j\bigr)
(q^{1+(j-i-1)\gamma}x_i/x_j)_{2\gamma-1} \\
\times 
\prod_{i=1}^{k_2} y_i^{\alpha_2}
(q^{1+(k_2-i)\gamma}y_i)_{\beta_2-1}
\prod_{1\leq i<j\leq k_2} y_j^{2\gamma}
\bigl(1-q^{(j-i)\gamma}y_i/y_j\bigr)
(q^{1+(j-i-1)\gamma}y_i/y_j)_{2\gamma-1} \\
\times
\prod_{i=1}^{k_1} \prod_{j=1}^{k_2} y_j^{-\gamma}
(q^{\beta_1+(k_1-k_2+j-i)\gamma}x_i/y_j)_{-\gamma} \\
=\prod_{i=1}^{k_1}\frac{\Gamma_q(\alpha_1+(k_1-i)\gamma+\la_i)
\Gamma_q(\beta_1+(i-k_2-1)\gamma)\Gamma_q(i\gamma)}
{\Gamma_q(\alpha_1+\beta_1+(2k_1-k_2-i-1)\gamma+\la_i)\Gamma_q(\gamma)} \\
\times \prod_{i=1}^{k_2}
\frac{\Gamma_q(\alpha_2+(k_2-i)\gamma+\mu_i)
\Gamma_q(\beta_2+(i-1)\gamma)\Gamma_q(i\gamma)}
{\Gamma_q(\alpha_2+\beta_2+(2k_2-k_1-i-1)\gamma+\mu_i)\Gamma_q(\gamma)} \\
\times\prod_{i=1}^{k_1} \prod_{j=1}^{k_2} 
\frac{\Gamma_q(\alpha_1+\alpha_2+(k_1+k_2-i-j-1)\gamma+\la_i+\mu_j)}
{\Gamma_q(\alpha_1+\alpha_2+(k_1+k_2-i-j)\gamma+\la_i+\mu_j)}.
\end{multline*}
Here $x_i:=q^{\eta_i}$ and $y_i:=q^{\nu_i}$, so that
\begin{equation}\label{xy}
0<x_1<\cdots<x_{k_1}\leq 1\and 0<y_1<\cdots<y_{k_1}\leq 1.
\end{equation}
The above is essentially a $(k_1+k_2)$-dimensional $q$-integral 
(more on this in the next section) and all that remains is
to let $q$ tend to $1$ from below.
The resulting integrand, however, depends sensitively on the relative 
ordering between the $x_i$ and $y_j$.
Indeed \cite{W09},
\begin{multline*}
\lim_{q\to 1^{-}}
y_j^{-\gamma} (q^{\beta_1+(k_1-k_2+j-i)\gamma}x_i/y_j)_{-\gamma}  \\
=\abs{x_i-y_j}^{-\gamma} \times
\begin{cases}1 & \text{if $x_i<y_j$} \\[2mm]
\displaystyle
\frac{\sin\pi (\beta_1-(i-j-k_1+k_2)\gamma)}
{\sin\pi (\beta_1-(i-j-k_1+k_2+1)\gamma)}
& \text{if $x_i>y_j$}.
\end{cases}
\end{multline*}
Consequently, before we can take the required limit we must fix
a complete ordering among the integration variables (compatible with
\eqref{xy}) and sum over all admissible orderings.
This is exactly what is done at the beginning of this section
and in the remainder we assume that 
\[
(x_1,\dots,x_{k_1},y_1,\dots,y_{k_2})\in I_{a_1,\dots,a_{k_1}}^{k_1,k_2}[0,1].
\]
To find how to weigh this domain we recall that according to 
\eqref{order} $y_{a_i+1}>x_i>y_{a_i}>\dots>y_1$.
The correct weight is thus
\[
\prod_{i=1}^{k_1} \prod_{j=1}^{a_i} 
\frac{\sin\pi (\beta_1-(i-j-k_1+k_2)\gamma)}
{\sin\pi (\beta_1-(i-j-k_1+k_2+1)\gamma)}
=\prod_{i=1}^{k_1} 
\frac{\sin\pi (\beta_1-(i-a_i-k_1+k_2)\gamma)}
{\sin\pi (\beta_1-(i-k_1+k_2)\gamma)},
\]
in accordance with $C^{k_1,k_2}_{\beta_1,\gamma}[0,1]$, see \eqref{Ckk}.
\end{proof}

\section{The Askey--Habsieger--Kadell integral}

For $0<q<1$ the $q$-integral on $[0,1]$ is defined as
\begin{equation}\label{qint}
\int_0^1 f(x)\dupq x=(1-q)\sum_{k=0}^{\infty} f(q^k) q^k,
\end{equation}
where it is assumed the series on the right converges.
When $q\to 1^-$ the $q$-integral reduces, at least formally,
to the Riemann integral of $f$ on the unit interval.
An obvious $n$-dimensional analogue of \eqref{qint} is
\[
\Int_{[0,1]^n} f(X)\dupq X=(1-q)^n\sum_{k_1,\dots,k_n=0}^{\infty} 
f(q^{k_1},\dots,q^{k_n}) q^{k_1+\cdots+k_n},
\]
where the multiple sum on the right is assumed to be absolutely convergent
and where $f(X)=f(x_1,\dots,x_n)$ and $\dupq X=\dupq x_1\cdots\dupq
x_n$.

In 1980 Askey \cite{Askey80} conjectured a $q$-analogue of the
Selberg integral when the parameter $\gamma$ is a nonnegative integer,
say $k$:
\begin{multline*}
\Int_{[0,1]^n}\prod_{i=1}^n x_i^{\alpha-1} 
(x_iq)_{\beta-1}
\prod_{1\leq i<j\leq n} x_i^{2k} (q^{1-k}x_j/x_i)_{2k}
\,\dupq X \\
=q^{\alpha k\binom{n}{2}+2k^2\binom{n}{3}}
\prod_{i=0}^{n-1} 
\frac{\Gamma_q(\alpha+ik)\Gamma_q(\beta+ik)\Gamma_q(1+(i+1)k)}
{\Gamma_q(\alpha+\beta+(n+i-1)k)\Gamma_q(1+k)} ,
\end{multline*}
for $\textup{Re}(\alpha)>0$ and $\beta\neq 0,-1,-2,\dots$. 
Askey's conjecture was proved independently by Habsieger \cite{Habsieger88}
and Kadell \cite{Kadell88}.

Just as the ordinary Selberg integral, the Askey--Habsieger--Kadell integral 
can be generalised by the inclusion of symmetric functions in the
integrand. Specifically, Kaneko \cite{Kaneko96} and Macdonald 
\cite{Macdonald95} proved that
\begin{multline}\label{qKM}
\Int_{[0,1]^n} \tilde{\P}_{\la}(X;q,q^k)
\prod_{i=1}^n x_i^{\alpha-1}(x_iq)_{\beta-1}
\prod_{1\leq i<j\leq n} x_i^{2k}
(q^{1-k}x_j/x_i)_{2k}\, \dupq X \\
=q^{\alpha k\binom{n}{2}+2k^2\binom{n}{3}}
\prod_{i=1}^n
\frac{\Gamma_q(\alpha+(n-i)k+\la_i)\Gamma_q(\beta+(i-1)k)\Gamma_q(ik+1)}
{\Gamma_q(\alpha+\beta+(2n-i-1)k+\la_i)\Gamma_q(k+1)},
\end{multline}
for $\textup{Re}(\alpha)>-\la_n$ and $\beta\neq 0,-1,-2,\dots$.

The $\sln{n}$--\,$\sln{m}$ transformation formula of 
Theorem~\ref{thmPhitrafo} allows the for the Askey--Habsieger--Kadell 
integral as well as its generalisation \eqref{qKM} to be extended to
a transformation between integrals of different dimensions.
For $\la=(\la_1,\dots,\la_n)$ and $\mu=(\mu_1,\dots,\mu_m)$ define
\begin{multline*}
S^{(n,m)}_{\la\mu}(\alpha_1,\alpha_2,\beta;k)
=\Int_{[0,1]^n}\tilde{\P}_{\la}(X;q,q^k)\prod_{i=1}^n x_i^{\alpha_1-1} 
(x_iq)_{\beta-(n-1)k-1}\\
\times \prod_{i=1}^n \prod_{j=1}^m 
\frac{(x_iq)_{\alpha_2+\beta+\mu_j+(m-n-j)k-1}}
{(x_iq)_{\alpha_2+\beta+\mu_j+(m-n-j+1)k-1}} 
\prod_{1\leq i<j\leq n} x_i^{2k} (q^{1-k}x_j/x_i)_{2k}
\,\textup{d}_q X.
\end{multline*}

\begin{theorem}
Let $\la=(\la_1,\dots,\la_n)$, $\mu=(\mu_1,\dots,\mu_m)$ be partitions,
$k$ a nonnegative integer and $\alpha_1,\alpha_2,\beta\in\Complex$.
Then
\begin{align*}
S^{(n,m)}_{\la\mu}(\alpha_1,\alpha_2,\beta;k)
&=q^{\alpha_1 k\binom{n}{2}-\alpha_2 k\binom{m}{2}+
2k^2\binom{n}{3}-2k^2\binom{m}{3}} \,
S^{(m,n)}_{\mu\la}(\alpha_2,\alpha_1,\beta;k) \\
&\quad \times \prod_{i=1}^n 
\frac{\Gamma_q(\beta-(i-1)k)\Gamma_q(\alpha_1+\la_i+(n-i)k)\Gamma_q(ik+1)}
{\Gamma_q(\alpha_1+\beta+\la_i+(n-m-i)k)\Gamma_q(k+1)} \\ 
&\quad \times \prod_{i=1}^m 
\frac{\Gamma_q(\alpha_2+\beta+\mu_i+(m-n-i)k)\Gamma_q(k+1)}
{\Gamma_q(\beta-(i-1)k)\Gamma_q(\alpha_2+\mu_i+(m-i)k)\Gamma_q(ik+1)}
\end{align*}
for $\textup{Re}(\alpha_1)>-\la_n$, $\textup{Re}(\alpha_2)>-\mu_m$, and
generic $\beta$.
\end{theorem}
By ``generic $\beta$'' it is meant that $\beta$ should avoid
a countable set of isolated singularities. More precisely, $\beta$ 
should be such that none of $\beta-(n-1)k$, $\beta-(m-1)k$, 
$\alpha_1+\beta+\la_j+(n-m-j)k$ and
$\alpha_2+\beta+\mu_j+(m-n-j)k$ take nonpositive integer values.
Since $S^{(0,n)}_{0,\la}(\alpha_2,\alpha_1,\beta;k)=1$ the $m=0$ case of the 
theorem corresponds to \eqref{qKM} with 
$(\alpha,\beta)\mapsto (\alpha_1,\beta-(n-1)k)$.

\begin{proof}
The method of proof is identical to that employed in 
\cite[Theorem 1.1]{W05} and we only sketch the details
of what are essentially elementary manipulations.

We specialise $X\mapsto c\spec{\la}_n$ and $Y\mapsto b\spec{\mu}_m$
in Theorem~\ref{thmPhitrafo} and apply the evaluation symmetry \eqref{TK}
to obtain
\begin{multline*}
\sum_{\nu} 
\frac{(a,abq^{\mu_1}t^{m-2},\dots,abq^{\mu_m}t^{-1})_{\nu}}
{(abq^{\mu_1}t^{m-1},\dots,abq^{\mu_m})_{\nu}}\,
\P_{\nu}(c\spec{0}_n)\tilde{\P}_{\la}(\spec{\nu}_n)\\
=\biggl(\:\prod_{i=1}^n \frac{(acq^{\la_i}t^{n-i})_{\infty}}
{(cq^{\la_i}t^{n-i})_{\infty}}
\biggr)\biggl(\:\prod_{i=1}^m \frac{(bq^{\mu_i}t^{m-i})_{\infty}}
{(abq^{\mu_i}t^{m-i})_{\infty}} \biggr) \\
\times \sum_{\nu} 
\frac{(a,acq^{\la_1}t^{n-2},\dots,acq^{\la_n}t^{-1})_{\nu}}
{(acq^{\la_1}t^{n-1},\dots,acq^{\la_n})_{\nu}}\,
\P_{\nu}(b\spec{0}_m)\tilde{\P}_{\mu}(\spec{\nu}_m).
\end{multline*}

Next we replace $t\mapsto q^k$ with $k$ a positive integer
and replace $a\mapsto q^{\beta}$, $b\mapsto q^{\alpha_2}$ and $c\mapsto
q^{\alpha_1}$. Then we apply 
\cite[Lemma 3.1]{W05} to write the $\nu$-sums as $n$-fold
unrestricted sums, and the claim follows.
\end{proof}

The above derivation can be repeated starting
from Theorem~\ref{thmCauchy}. The result is an $\sln{3}$ variant of the
$q$-integral \eqref{qKM}. Problem with the theorem below is, however, that
it does not converge in the $q\to 1^{-}$ limit unless $m$ or $n$ is $0$. 
(This can be remedied by
replacing $[0,1]^{m+n}$ by appropriate multiple Pochhammer double loops).
\begin{theorem}
Let $\la=(\la_1,\dots,\la_n)$, $\mu=(\mu_1,\dots,\mu_m)$ be partitions,
$k$ a nonnegative integer and $\alpha_1,\alpha_2,\beta_1,\beta_2\in\Complex$ 
such that $\beta_1+\beta_2=k+1$. Then
\begin{align*}
\Int_{[0,1]^{n+m}} & \tilde{\P}_{\la}(X;q,q^k) 
\prod_{i=1}^n x_i^{\alpha_1} (qx_i)_{\beta_1-1}
\prod_{1\leq i<j\leq n} x_j^{2k} (q^{1-k} x_i/x_j)_{2k} \\[-2mm]
\times\, & \tilde{\P}_{\mu}(Y;q,q^k)
\prod_{i=1}^m y_i^{\alpha_2} (qy_i)_{\beta_2-1}
\prod_{1\leq i<j\leq m} y_j^{2k} (q^{1-k} y_i/y_j)_{2k} \\
&\hspace{4cm}\times\prod_{i=1}^n \prod_{j=1}^m 
y_j^{-k} (q^{\beta_1}x_i/y_j)_{-k} 
\: \dupq X \, \dupq Y \\
&=q^{\alpha_1 k\binom{n}{2}+\alpha_2 k\binom{m}{2}+
2k^2\binom{n}{3}+2k^2\binom{n}{3}-k^2 n\binom{m}{2}} \\
&\quad\times
\prod_{i=1}^n\frac{\Gamma_q(\alpha_1+(n-i)k+\la_i)
\Gamma_q(\beta_1+(i-m-1)k)\Gamma_q(ik+1)}
{\Gamma_q(\alpha_1+\beta_1+(2n-m-i-1)k+\la_i)\Gamma_q(k)} \\
&\quad\times
\prod_{i=1}^m \frac{\Gamma_q(\alpha_2+(m-i)k+\mu_i)
\Gamma_q(\beta_2+(i-1)k)\Gamma_q(ik+1)}
{\Gamma_q(\alpha_2+\beta_2+(2m-n-i-1)k+\mu_i)\Gamma_q(k)} \\
&\quad\times\prod_{i=1}^n \prod_{j=1}^m
\frac{\Gamma_q(\alpha_1+\alpha_2+(n+m-i-j-1)k+\la_i+\mu_j)}
{\Gamma_q(\alpha_1+\alpha_2+(n+m-i-j)k+\la_i+\mu_j)},
\end{align*}
for $\textup{Re}(\alpha_1)>-\la_n$,
for $\textup{Re}(\alpha_2)>-\mu_m$,
and $\beta_1,\beta_2\neq 0,-1,-2,\dots$.
\end{theorem}

\subsection*{Acknowledgements}
I thank Eric Rains for helpful discussions.

\bibliographystyle{amsplain}

\end{document}